\documentclass{article}

\usepackage{arxiv}

\usepackage[utf8]{inputenc} 
\usepackage[T1]{fontenc}    
\usepackage{hyperref}       
\usepackage{url}            
\usepackage{booktabs}       
\usepackage{amsfonts}       
\usepackage{nicefrac}       
\usepackage{microtype}      
\usepackage{lipsum}
\usepackage{amsthm}
\usepackage[fleqn]{amsmath}
\usepackage[ruled]{algorithm2e}
\usepackage{multirow}
\usepackage{subfigure}
\usepackage{epsfig}
\newtheorem{Proposition}{Proposition}

\title{Last-mile Delivery: Optimal Locker Location Under Multinomial Logit Choice Model}

\author{
  Yun Hui Lin \\
  National University of Singapore\\
  \texttt{isemlyh@gmail.com} \\
   \And
 Dongdong He\\
  National University of Singapore\\
  \texttt{isehd@nus.edu.sg} \\
     \And
 Yuan Wang\\
  National University of Singapore\\
  \texttt{iseway@nus.edu.sg} \\
      \And
 Loo Hay~Lee\\
  National University of Singapore\\
  \texttt{iseleelh@nus.edu.sg} \\ 
}

\begin{document}
\maketitle

\begin{abstract}
One innovative solution to the last-mile delivery problem is the self-service locker system. Motivated by a real case in Singapore, we consider a POP-Locker Alliance who operates a set of POP-stations and wishes to improve the last-mile delivery by opening new locker facilities. We propose a quantitative approach to determine the optimal locker location with the objective to maximize the overall service provided by the alliance. Customer's choices regarding the use of facilities are explicitly considered. They are predicted by a multinomial logit model. We then formulate the location problem as a multi-ratio linear-fractional 0-1 program and provide two solution approaches. The first one is to reformulate the original problem as a mixed-integer linear program, which is further strengthened using conditional McCormick inequalities. This approach is an exact method, developed for small-scale problems. For large-scale problems, we propose a Suggest-and-Improve framework with two embedded algorithms. Numerical studies indicated that our framework is an efficient approach that yields high-quality solutions. Finally, we conducted a case study.  The results highlighted the importance of considering the customers’ choices. Under different parameter values of the multinomial logit model, the decisions could be completely different.  Therefore, the parameter value should be carefully estimated in advance. 
\end{abstract}


\section{Introduction}
\label{S:1}

In recent years, the e-commerce market has been expanding at a remarkable rate. E-commerce revenues worldwide are projected to grow to 6.54 trillion US dollars in 2022, according to Statista, a global business data platform \footnote{https://www.statista.com/statistics/379046/worldwide-retail-e-commerce-sales/}. The rapid growth of e-commerce has led to a sharp increase in the delivery of online purchasing parcels. The process of transporting goods to customers consists of many stages that have undergone significant improvements towards a more efficient and cost-effective manner. However, the final stage, i.e., delivering the parcels (goods) to the doorstep of the customer, is still \textit{not} efficient. In practice, the last-mile delivery (LMD) often involves one package per door (which disallows the economies of scale) and the couriers may have difficulties in finding the exact home address of the related customer. Meanwhile, the problem of ``Not-at-Home"  can result in an empty trip, especially when the delivery needs a signature of a receipt to confirm the success of a delivery~\cite{Deutsch2018}. \cite{Kedia2017} showed that missing a delivery may also require a pick-up from a collection point with unsuitable location or time requirements. These issues lower the efficiency of the LMD and make the delivery process costly. Therefore, how to address these issues is an essential topic for e-commerce delivery. As a response, operators have started to offer alternative options to regular home deliveries.

Self-service technology (SST) provides an innovative solution. In recent years, self-service POP-stations and parcel lockers 
have been deployed for the LMD. They allow logistics companies to deliver parcels to the facilities with economies of scale. Customers can collect the parcels at their convenience.  Recent research showed that self-service POP-stations and parcel locker networks are gaining more popularities~\cite{Dieke2013,Morganti2014} and their deployments significantly improve the consumers'  shopping experience and enhance competitiveness and performance of the business~ \cite{Vakulenko2018}. In Singapore, Singapore Post has launched the ``Locker Alliance" program. The alliance provides an extensive delivery network that consists of self-service POP-stations and parcel lockers. Using this system, logistics companies can avoid the ``Not-at-Home" problem in deliveries and customers do not need to wait at home for receiving the parcels. Therefore, it creates a ``win-win" situation for companies and customers.


Recently, there is a proliferation of research on the parcel locker. Many studies have highlighted the benefits of using it from both economic and environmental aspects.  For instance, \cite{Kaemaeraeinen2001} used delivery data from the suburban area of Helsinki and compared the delivery costs of two strategies, namely, the regular home delivery and the batch delivery to reception boxes. Their results showed a 42\% cost reduction when the deliveries are made to boxes. Similarly, \cite{Punakivi2001} conducted a case study on a retail company in Finland and showed that the operational cost of regular home delivery is more expensive than that of delivery to reception boxes and parcel lockers.  \cite{Song2013} and \cite{Lemke2016} argued that the greenhouse gas emissions the LMD can be reduced by deploying the self-collection service (e.g., parcel lockers or collection-and-delivery points). Researchers have also examined the factors influencing customers' intention to use self-collection services. \cite{xu2013impact} showed that implementing customer segmentation in self-collection services can better suit customers' need and customers' satisfaction of home delivery service is not a significant factor that influences customers’ willingness to use self-collection services. Using survey data in Singapore, \cite{yuen2018investigation} argued that, to improve customers' intention, it is important is to integrate self-collection services into customers’ lifestyle, values and needs. Other studies focused on the conceptualization of hyperconnected lockers network~\cite{Faugere2017} and configuration designs of locker banks based on optimization models~\cite{Faugere2017a,Faugere2018}.  Among the literature, \cite{Deutsch2018} proposed a quantitative approach to determine the optimal locker location and the sizes of lockers. The paper assumed that customers made choices in a deterministic fashion and formulated an optimization problem to maximize the total profit of the operator.

In this paper, we study the optimal locker locations as well. Different from \cite{Deutsch2018}, we use a random utility model to predict customer's choices. In practice, customers may prefer some facility (not necessarily the most convenient one) according to their preferences. The preferences are, however, unknown to the operators~\cite{ljubic2018outer}. Empirical and experimental studies also show that individual choice behaviors are difficult to capture and they are typically interpreted probabilistically~\cite{gul2014random}. One way to predict the customer's choices is to apply the random utility theory.  More specifically, we use the multinomial logit (MNL) modell because we have witnessed its successful applications in various facility location problems where the customer's choices are considered~\cite{Haase2014comparison,ljubic2018outer}. Besides, motivated by the ``Locker Alliance" program in Singapore, we consider an alliance scenario where two operators corporate with each other to maximize the overall service level. The problem is not only to choose the candidate locations where the new lockers will be built, but also to close some unnecessary facilities whose deployments can reduce the overall service level when new lockers come to the network.

Our proposed model is a multi-ratio linear-fractional 0-1 program (MLFP).  Indeed, the MLFP has been applied to many economic and engineering optimization problems, including retail assortment planning~\cite{mendez2014branch}, set covering~\cite{amaldi2012hyperbolic}, design of stochastic service systems~\cite{elhedhli2005exact}, and competitive facility location problems~\cite{ljubic2018outer}. To solve the problem, the most straightforward approach is to reformulate it as an equivalent mixed-integer linear programming (MILP)~\cite{Haase2014comparison,mendez2014branch}. 
The major drawback of this approach is that it cannot scale well on problem size,which may fail to solve a large-scale (even just a moderate-scale) problem in a reasonable time. On the other hand, one exact approach that is capable of solving large-scale MLFPs is the outer approximation algorithm (OA)~\cite{elhedhli2005exact,ljubic2018outer}, which proceeds by solving a sequence of sparse MILPs with an increasing number of constraints. The OA algorithm guarantees finite convergence and outperforms the MILP approach by a large margin in general. However, it requires the concavity of the objective function, which makes it not valid for our problem because our objective function is not concave. Therefore, the quadratic transform approach, which transforms the MLFP into a bi-concave maximization program (BCMP), has been developed to solve such issue~\cite{benson2004global1,shen2018fractional}. But solving the BCMP itself requires significant efforts to design a specific algorithm~\cite{benson2007solving} and the solving process can be time-consuming as well.
For comprehensive survey papers on the applications and algorithms for MLFPs, we refer the reader to \cite{schaible2003fractional} and \cite{borrero2017fractional}. In this paper, we apply the quadratic transform approach and design an efficient ``Suggest-and-Improve" framework with two embedded algorithms to solve the large-scale BCMP.

To summarize, our contributions are fourfold. The first one is that we consider a POP-Locker Alliance who operates a set of POP-stations and wishes to improve the service level by opening new locker facilities and closing unnecessary (if any) POP-stations. We predict customer's choices by a MNL model and propose a MLFP to determine the optimal locker location. To the best of our knowledge, this paper is the first to incorporate the random utility theory into the locker location problem.  The proposed model is a non-concave maximization model. Our second contribution is to reformulate the original problem as an equivalent MILP, which is further strengthened using McCormick inequalities. This approach is an exact method that works well for small-scale problems. However, it cannot scale well on the problem size.  The poor performance of MILPs on large-scale problems motivates our third contribution, i.e., the ``Suggest-and-Improve" framework. In the Suggest procedure, we propose the \textit{Quadratic Transform with Linear Alternating} (QT-LA) algorithm to find a high-quality solution. The solution by QT-LA is subsequently improved by a searching algorithm in Improve procedure.  Numerical studies indicated that our proposed framework is an efficient approach that yields high-quality solutions. Finally, we conduct a case study based on the POP-locker system in Singapore.

The rest of the paper is organized as follows. Section~\ref{s:pd} covers the problem description and the basic mathematical model. We present the MILP solution method in Section~\ref{s:milp} and the Suggest-and-Improve framework in Section~\ref{s:sif}. Numerical studies are conducted in Section~\ref{s:ns}. We conclude the paper in Section~\ref{s:conclusion}.

\section{Problem description}
\label{s:pd}

In this section, we present the problem description. Table~\ref{tab:Nomenclature} gives the notations. Additional notations will be introduced when necessary.

\begin{table} \label{tab:Nomenclature}
\caption{Nomenclature}
\centering
\label{tab:Nomenclature}       
\begin{tabular}{lll}
\hline\noalign{\smallskip}
\multicolumn{3}{l}{\underline{Sets}}\\
$I$ &      $:$ & set of customer zones\\
$J$ &     $:$ & set of candidate locker facilities\\
$K$ &     $:$ & set of POP-stations\\

\multicolumn{3}{l}{\underline{Parameters}}\\
$l_{ik}$ &     $:$ & distance between demand zone $i$ and pos-station $k$, $\forall i \in I, k \in K$\\
$L_{ij}$ &     $:$ & distance between demand zone $i$ and locker $j$, $\forall i \in I, j \in J$\\
$\tilde{a}_{ik}$ &     $:$ & service level when a customer in zone $i$ uses POP-station $k$. $\tilde{a}_{ik} \in [0,1]$,$\forall i \in I, k \in K$\\
$a_{ij}$ &     $:$ & service level when a customer in zone $i$ uses locker $j$. $a_{ij} \in [0,1]$, $\forall i \in I, j \in J$\\
$v_{im}$ &     $:$ & utility when a customer in zone $i$ uses facility $m$, $\forall m \in J \cup K$\\
$\tilde{\theta}_{ik}$ &     $:$ & preference a customer in zone $i$ to POP-station $k$, $\forall k \in K$\\
$\theta_{ij}$ &     $:$ & preference a customer in zone $i$ to locker $j$, $\forall j \in J$\\
$d_i$ &     $:$ &proportion of the total demand that arises from zone $i$, $\forall i \in I$ \\
$P$ &     $:$ &maximum number of open lockers/ maximum number of closed POP-stations \\

\multicolumn{3}{l}{\underline{Decision Variables}}\\
$x_j$ & $:$ & 1, if locker $j$ is open; 0, otherwise. $\forall i \in I$\\
$r_k$ & $:$ & 1, if POP-station $k$ remains open; 0, otherwise. $\forall k \in K$\\
\noalign{\smallskip}\hline
\end{tabular}
\end{table}

Consider a POP-Locker Alliance that operates a set of POP-stations (denoted by set $K$) to serve a set of customer zones (denoted by set $I$).  The customers are assumed homogeneous in their observable characteristics. The distance between customer zone $i$ and POP-station $k$ is $l_{ik}$, $\forall i \in I, k \in K$. The service level that a customer in zone $i$ patronizes POP-station $k$ is measured by $\tilde{a}_{ik} \in [0,1]$. We assume $\tilde{a}_{ik}$ is a nondecreasing function of $l_{ik}$, i.e., when the distance between node $i$ and node $k$ increases, the service level should not increase. Besides POP-stations, there are potential locker facilities (denoted by set $J$). The distance between customer zone $i$ and locker $j$ is $L_{ij}$, $\forall i \in I, j \in J$. The service level that a customer in zone $i$ patronizes locker $j$ is measured by $a_{ij} \in [0,1]$. Similarly, we assume $a_{ij}$ is a nondecreasing function of $L_{ij}$.

Now, the alliance wishes to open new locker facilities and to close some unnecessary (if any) POP-stations so that the overall service level provided by the POP-Locker system is maximized. If the alliance assigns customers to facilities or customers always select the nearest facility, the problem becomes a type of classical coverage problems~\cite{daskin2011network}. It is then advantageous to open P locker facilities and not to close any POP-station when we restrict the maximum number of open locker facilities to P. However, customers are taking the initiative in e-commence. They may prefer a facility that is not the nearest according to their preferences. The preferences are, however, unknown and hard to observe~\cite{ljubic2018outer}. Empirical and experimental evidences have supported that individual choice behaviors should be interpreted probabilistically~\cite{gul2014random}. Therefore, the results from the traditional coverage model may not be relevant.

The multinomial logit model (MNL) provides a way to forecast customers' discrete choice behaviors. It assumes that customers in zone $i$ will maximize their utility:
\begin{align}
u_{im} = v_{im} + \epsilon_{im},~\forall i \in I, m \in J\cup K
\end{align}
when choosing a facility. Here, $v_{im}$ is the deterministic part of the utility, which is measured as a function of convenience (distance). $\epsilon_{im}$ is a random term that is independent identically extreme value distributed.

Define a binary variable $r_k$, $\forall k \in K$, such that $r_k = 1$, if POP-station  $k$ remains in operation; $r_k = 0$, otherwise. Define a binary variable $x_j$, $\forall j \in J$, such that $x_j = 1$, if locker $j$ is open; $x_j = 0$, otherwise. According to MNL, the probability that a customer in zone $i$ selects POP-station $k$ is
\begin{align}
p_{ik} = \frac{e^{v_{ik}}r_k}{\sum_{k \in K}e^{v_{ik}}r_k +  \sum_{j \in J} e^{v_{ij}} x_j }
\end{align}
Similarly, the probability that a customer in zone $i$ selects locker $j$ is
\begin{align}
p_{ij} = \frac{e^{v_{ij}} x_j}{\sum_{k \in K}e^{v_{ik}}r_k +  \sum_{j \in J} e^{v_{ij}} x_j }
\end{align}
Let $\tilde{\theta}_{ik} = e^{v_{ik}}$,$\forall k \in K$ and  $\theta_{ij} = e^{v_{ij}}$,$\forall i \in I$. The number of customers in $i$ that use POP-station $k$ is:
 \begin{align}
D_{ik} = d_j p_{ik} = \frac{d_i \tilde{\theta}_{ik}r_k}{\sum_{k \in K}\tilde{\theta}_{ik}r_k +  \sum_{j \in J} \theta_{ij} x_j }
\end{align}
Similarly, the number of customers in $i$ that use locker $j$ is
\begin{align}
D_{ij} = d_j p_{ij} = \frac{d_i \theta_{ij} x_j}{\sum_{k \in K}\tilde{\theta}_{ik}r_k +  \sum_{j \in J} \theta_{ij} x_j }
\end{align}
We then compute the overall service level by
\begin{align}
\label{eqt:obj} C(x,r) =\sum_{i \in I}(\sum_{k \in K}\tilde{a}_{ik}D_{ik} +\sum_{j \in J}a_{ij}D_{ij}) = \sum_{i \in I} \frac{d_i (\sum_{k \in K} \tilde{a}_{ik}\tilde{\theta}_{ik}r_k + \sum_{j \in J} a_{ij} \theta_{ij}x_j)}{\sum_{k \in K}\tilde{\theta}_{ik}r_k +  \sum_{j \in J} \theta_{ij} x_j }
\end{align}
It is easy to see that $C(x,r)$ lies between $0$ and $1$. To facilitate later discussions, we define parameters $\tilde{b}_{ik}$ and $b_{ij}$ such that $\tilde{b}_{ik} = d_i \tilde{a}_{ik} \tilde{\theta}_{ik}$ and $b_{ij} = d_i a_{ij} \theta_{ij}$. With them, we propose the following model for determining the optimal locker location under MNL model:
\begin{align}
\label{eqt:p1_obj}\textit{\textbf{[P1]}} \quad  \max~& C(x,r) = \sum_{i \in I} \frac{ \sum_{k \in K} \tilde{b}_{ik}r_k + \sum_{j \in J} b_{ij}x_j}{\sum_{k \in K}\tilde{\theta}_{ik}r_k +  \sum_{j \in J} \theta_{ij} x_j }\\
 \label{constr:P1-1} st.~& \sum_{j \in J} x_j \leq P \\
 \label{constr:P1-2}  & \sum_{k \in K}(1- r_k) \leq P \\
 \label{constr:P1-3} & x_j \in \{0,1\}, \forall j \in J \\
  \label{constr:P1-4} & r_k \in \{0,1\}, \forall k \in K
\end{align}
where $0 \leq P \leq \max \{|J|,|K|\}$. The objective is to maximize the overall service level $C(x,r)$. Constraint~(\ref{constr:P1-1}) imposes that the maximum number of open lockers is $P$. Constraint~(\ref{constr:P1-2}) indicates that the maximum number of closed POP-stations cannot exceed $P$. In practice, the maximum number of open lockers and  the maximum number of closed POP-stations need not be the same. For simplicity and without loss of generality, we assume a same value for both constraints.  Meanwhile, if constraint~(\ref{constr:P1-1}) and (\ref{constr:P1-2}) are set to equalities, then they simply indicate $P$ POP-stations will be replaced by the same number of lockers in the network. Finally, constraints~(\ref{constr:P1-3}) and~(\ref{constr:P1-4}) specify that $x$ and $r$ are binary variables.

The proposed model is a constrained MLFP. It is a non-concave 0-1 maximization problem in general. Such a problem is NP-hard if $|I|\geq 2$ ~\cite{prokopyev2005multiple}. In the next two sections, we will provide both exact and heuristic approaches.

\section{Mixed-integer linear programming approach}
\label{s:milp}

In this section, we present an exact solution approach that works well for small-scale problems. We first derive an equivalent MILP formulation  for P1 and then strengthen the formulation by conditional McCormick inequalities.

\subsection{Basic formulation}

The MILP reformulation approach is widely used for multi-ratio linear fractional 0-1 programing~\cite{bront2009column,Haase2014comparison,mendez2014branch}. Given recent improvements of off-the-shelf MILP solvers such as CPLEX and Gurobi, reformulating P1 as its equivalent MILP provides a standard solution framework that is straightforward to implement.

Let $z_i = 1/(\sum_{k \in K}\tilde{\theta}_{ik}r_k +  \sum_{j \in J} \theta_{ij} x_j)$. P1 is equivalent to the bilinear mixed 0-1 problem:
\begin{align}
 \qquad   \qquad  \max~& \sum_{i \in I} \sum_{k \in K} \tilde{b}_{ik} r_k z_i + \sum_{i \in I} \sum_{j \in J} b_{ij}x_jz_i\\
 \label{constr:P2-1} st.~&\sum_{k \in K}\tilde{\theta}_{ik}r_k z_i +  \sum_{j \in J} \theta_{ij} x_jz_i = 1,~\forall i \in I \\
 \label{constr:P2-2} & z_i \geq 0,~\forall i \in I  \\
\nonumber  & (\ref{constr:P1-1}) - (\ref{constr:P1-4})
\end{align}
which can be recast as a MILP by linearizing the bilinear terms $x_j z_i$ and $r_k z_i$. Let $y_{ij} = x_j z_i$ and $Y_{ik} = r_k z_i$. Consider the following MILP:
\begin{align}
\textit{\textbf{[MILP]}}\quad \max~& \sum_{i \in I}\sum_{k \in K} \tilde{b}_{ik}Y_{ik} + \sum_{i \in I}\sum_{j \in J} b_{ij}y_{ij}\\
 \label{constr:MILP-1} st.~&\sum_{k \in K}\tilde{\theta}_{ik}Y_{ik}  +  \sum_{j \in J} \theta_{ij} y_{ij} = 1,~\forall i \in I \\
 \label{constr:MILP-2} & z_i - U(1-x_j) \leq y_{ij}  \leq z_i,~\forall i \in I, j \in J \\
 \label{constr:MILP-3} &0 \leq y_{ij} \leq U x_j,~\forall i \in I, j \in J \\
 \label{constr:MILP-4} & z_i - U(1-r_k) \leq Y_{ik}  \leq z_i,~\forall i \in I, k \in K \\
 \label{constr:MILP-5} &0 \leq Y_{ik} \leq U r_k,~\forall i \in I, k \in K \\
 \nonumber  & (\ref{constr:P1-1}) - (\ref{constr:P1-4}), (\ref{constr:P2-2})
\end{align}
where $U$ be a sufficiently large number. It can easily verified that this MILP is equivalent to P1. Specifically, (\ref{constr:MILP-2})-(\ref{constr:MILP-3}) impose that $y_{ij} = x_jz_i$. When $x_j = 1$, $z_i = y_{ij}$ by (\ref{constr:MILP-2}); when $x_j = 0$, $z_i = 0$ by (\ref{constr:MILP-3}). Similarly, (\ref{constr:MILP-4})-(\ref{constr:MILP-5}) impose that $Y_{ik} = r_kz_i$.

This approach is straightforward. In theory, the MILP can yield the same optimal solution as P1. However, as shown in~\cite{bront2009column} and \cite{mendez2014branch}, the MILP does not scale well on the problem size. In practice, it may fail to yield the optimal solutions in reasonable computational times, even for moderately sized problems.

\subsection{Strengthened formulation: conditional McCormick inequalities}

To enhance the performance of the MILP, we strengthen the MILP formulation by using (conditional) McCormick estimators for the bilinear terms. We first look at the following proposition.
\begin{Proposition} \label{prop:mc}
A bilinear mixed 0-1 term $y = xz$, where $x$ is a 0-1 variable and $z$ is a continuous variable with $z > 0$, can be represented by the following linear inequalities: (a) $z - z^u \cdot (1- x) \leq y \leq z - z^l|_{x=0} \cdot (1-x)$; (b) $z^l|_{x=1} \cdot x \leq y \leq z^u|_{x=1} \cdot x$,
where $z^u$ is the global upper bound on $z$; $z^l|_{x=0}$ is the lower bound on $z$ when $x=0$; $z^l|_{x=1}$ and $z^u|_{x=1}$ is the lower bound and upper bound on $z$ when $x=1$.
\end{Proposition}
\begin{proof}
We need to check two cases:
\begin{itemize}
\item If $x=0$, we have: (a) $z - z^u \leq y \leq z - z^l|_{x=0}$; (b) $0 \leq y \leq 0$. Clearly, (a) is inactive because $- z^u \leq -z^l|_{x=0}$ and (b) enforces that $y=0$.
\item If $x=1$, we have: (a) $z \leq y \leq z $; (b) $z^l|_{x=1} \leq y \leq z^u|_{x=1}$. Clearly, (b) is inactive because $z^l|_{x=1} \leq z^u|_{x=1}$ and (a) enforces that $y=z$.
\end{itemize}
Combine the cases, inequalities (a) and (b) give $y = xz$.
\end{proof}
By Proposition~(\ref{prop:mc}), the bilinear terms $x_j z_i$ and $r_k z_i$ can be represented by linear inequalities with global and conditional bounds on $z_i$ subject to (\ref{constr:P1-1}) - (\ref{constr:P1-4}). 

We first give the global upper bound on $z_i$.
\begin{Proposition}
The following global upper bounds on $z_i$, $\forall i \in I$, is valid:
\begin{align}\label{eqt:global_upper}
z^u_i &= \frac{1}{\sum^{|K| - P}_{e=1}\tilde{\theta}_{i(e)}}
\end{align}
where $\tilde{\theta}_{i(e)}$ is the $e$th smallest value of $\tilde{\theta}_{in}$, $\forall n \in K$.
\end{Proposition}
\begin{proof}
See~\ref{appendix}
\end{proof}

Next, we derive the conditional bounds. Let $z^u_{i|x_j = \xi}$ and $z^l_{i|x_j = \xi}$ be the upper and lower bound on $z_i$ when an additional constraint $x_j = \xi$ is imposed.
\begin{Proposition}\label{prop:cond_x}
The following conditional bounds on $z_i$ given $x_j$ , $\forall i \in I, j \in J$, are valid:
\begin{align}
z^u_{i \mid x_j=1} &= \frac{1}{\sum^{|K| - P}_{e=1}\tilde{\theta}_{i(e)} +  \theta_{ij}} \\
z^l_{i \mid x_j=0} &= \frac{1}{\sum_{k \in K}\tilde{\theta}_{ik} +  \sum^{P}_{e=1} \theta_{i[e]}} \\
z^l_{i \mid x_j=1} &= \frac{1}{\sum_{k \in K}\tilde{\theta}_{ik} + \theta_{ij} + \sum^{P-1}_{e=1} \theta_{i[e]}}
\end{align}
where $\tilde{\theta}_{i(e)}$ is the $e$th smallest value of $\tilde{\theta}_{in}$, $\forall n \in K$; $\theta_{i[e]}$ is the $e$th largest value of $\theta_{im}$, $\forall m \in J \backslash \{j\}$.
\end{Proposition}
\begin{proof}
See~\ref{appendix}
\end{proof}
Similarly, let $z^u_{i|r_k = \xi}$ and $z^l_{i|r_k = \xi}$ be the upper and lower bound on $z_i$ when an additional constraint $r_k = \xi$ is imposed.

\begin{Proposition}
The following conditional bounds on $z_i$ given $r_k$ , $\forall i \in I, k \in K$, are valid:
\begin{align}
z^u_{i \mid r_k=1} &= \frac{1}{\sum^{|K| - P -1}_{e=1}\tilde{\theta}_{i(e)} +  \tilde{\theta}_{ik}} \\
z^l_{i \mid r_k=0} &= \frac{1}{\sum_{k \in K\backslash \{k\}}\tilde{\theta}_{ik} +  \sum^{P}_{e=1} \theta_{i[e]}} \\
z^l_{i \mid r_k=1} &= \frac{1}{\sum_{k \in K}\tilde{\theta}_{ik} + \sum^{P}_{e=1} \theta_{i[e]}}
\end{align}
where $\tilde{\theta}_{i(e)}$ is the $e$th smallest value of $\tilde{\theta}_{in}$, $\forall n \in K$;  $\theta_{i[e]}$ is the $e$th largest value of $\theta_{im}$, $\forall m \in J$.
\end{Proposition}
\begin{proof}
See~\ref{appendix}
\end{proof}

With the above bounds on $z_i$, we can define the following (conditional) McCormick inequalities~\cite{mccormick1976} for $y_{ij} = x_jz_i$ and $Y_{ik} = r_k z_i$:
\begin{align}
\label{MC-1}                                 & z_i - z^u_i \cdot (1-x_j)  \leq  y_{ij}  \leq z_i - z^l_{i \mid x_j=0}\cdot (1-x_j),~\forall i \in I, j \in J \\
\label{MC-2} \textit{\textbf{[MC]}} \quad & z^l_{i \mid x_j=1} \cdot x_j \leq   y_{ij}  \leq z^u_{i \mid x_j=1} \cdot x_j,~\forall i \in I, j \in J \\
\label{MC-3}                                 & z_i - z^u_i \cdot (1-r_k) \leq Y_{ij} \leq z_i - z^l_{i \mid r_k=0}\cdot (1-r_k),~\forall i \in I, k \in K \\
\label{MC-4}                                 & z^l_{i \mid r_k=1} \cdot r_k  \leq Y_{ij} \leq z^u_{i \mid r_k=1}\cdot r_k,~\forall i \in I, k \in K
\end{align}
The idea of adding the MC is to tighten the continuous relaxation bound, thereby avoiding extensive branching and speeding up the computation~\cite{csen2018conic}. The impacts will be discussed in Section~\ref{s:numeric_mc}.

\section{Suggest-and-Improve framework}
\label{s:sif}

In Section~\ref{s:milp}, we present the MILP approach to P1. It is a favored approach for small-scale problems because it is an exact method that is easy to implement. However, the MILP has significantly more variables and constraints than P1. The solution times for large-scale problems can be prohibitive.  To solve large-scale problems, we introduce the ``Suggest-and-Improve" framework:
\begin{description}
\item[Step 1] \textit{Suggest Procedure}. Find a high-quality solution $(x,r)$ with a fast algorithm.
\item[Step 2] \textit{Improve Procedure}. Run a searching algorithm, using $(x,r)$ as the initial solution, to find a solution that is no worse than $(x,r)$.
\end{description}

We quote the name, ``Suggest-and-Improve", from~\cite{park2017general} who used this simple but flexible idea for nonconvex quadratically constrained quadratic programming. In this section, we will present methodologies for both procedures.

\subsection{Suggest procedure: QT-LA algorithm}


In Suggest procedure, we propose an algorithm that efficiently generates high-quality solutions. We start with the following proposition:
\begin{Proposition}
P1 is equivalent to the following problem:
\begin{align}
\label{eqt:qt_obj} \textbf{[P2]} \quad  \max~& \sum_{i \in I} 2 y_i \sqrt[]{\sum_{k \in K} \tilde{b}_{ik}r_k + \sum_{j \in J} b_{ij}x_j} - \sum_{i \in I}y_i^2( \sum_{k \in K} \tilde{\theta}_{ik}r_k +  \sum_{j \in J} \theta_{ij} x_j ) \\
 st.&~y_i \geq 0,~\forall i \in I \\
\nonumber  &(\ref{constr:P1-1}) - (\ref{constr:P1-4})
\end{align}
\end{Proposition}
\begin{proof}
See \cite{benson2004global1}.
\end{proof}
In effect, the equivalence of P1 and P2 can be easily established: we first optimize P2 over $y_i$. Clearly, this partial optimization is a concave maximization problem. By the first-order condition, we have:
\begin{align} \label{eqt:opt_y}
&y_i  := \frac{\sqrt[]{\sum_{k \in K} \tilde{b}_{ik}r_k + \sum_{j \in J} b_{ij}x_j}}{\sum_{k \in K}\tilde{\theta}_{ik}r_k +  \sum_{j \in J} \theta_{ij} x_j} , ~\forall i \in I
\end{align}
Plugging (\ref{eqt:opt_y}) into (\ref{eqt:qt_obj}) leads us to (\ref{eqt:p1_obj}).

According to~\cite{shen2018fractional}, P2 is also referred to as the Quadratic Transform (QT). Obviously, P2 is not a convex optimization problem. It is indeed a BCMP: when $x$ and $r$ are held fixed, the optimal $y$ can be found in closed form as (\ref{eqt:opt_y}). When $y$ is held fixed, P2 reduces to a (nonlinear) concave maximization program with variables $x$ and $r$.

One heuristic solution approach that leverages the bi-concavity of QT is to alternate between updating $y$ and solving the nonlinear program until some predetermined stopping condition is met, as proposed by~\cite{shen2018fractional}. However, this approach is initially proposed for continuous problems. When it is applied to 0-1 problems, it may stop at low-quality solutions or even loop around candidate solutions without convergence. One reason for its poor performance is the ``aggressive" rule for updating $y$: $y$ only depends on by the current solution of $x$ and $r$  by (\ref{eqt:opt_y}). Meanwhile, in our problem, solving the nonlinear 0-1 program for multiple iterations can be expensive when the problem size is large.

To overcome the issues, we propose a new alternating algorithm. Define function $f_i(x,r)$ such that
\begin{align}
f_i(x,r) = \sqrt[]{\sum_{k \in K} \tilde{b}_{ik}r_k + \sum_{j \in J} b_{ij}x_j}
\end{align}
We can rewrite P2 in the hypograph form of $f_i$, i.e.,
\begin{align}
\textit{\textbf{[P3]}} \quad  \max~& \sum_{i \in I} 2 y_i \beta_i - \sum_{i \in I}y_i^2( \sum_{k \in K} \tilde{\theta}_{ik}r_k +  \sum_{j \in J} \theta_{ij} x_j ) \\
\label{eqt:beta} st.~~&\beta_i \leq f_i(x,r) \\
 \nonumber &(\ref{constr:P1-1}) - (\ref{constr:P1-4})
\end{align}
Now, given any point $(\bar{x},\bar{r})$, since $f_i(x,r)$ is a concave function, we can bound it from above by its first-order linear approximation on $(\bar{x},\bar{r})$. The following constraint is thus valid for P3:
\begin{align}
\label{eqt:valid} \beta_i \leq  \sum_{j \in J} \frac{\partial f_i(\bar{x},\bar{r})}{\partial x_j}(x_j - \bar{x}_j) +   \sum_{k \in K} \frac{\partial f_i(\bar{x},\bar{r})}{\partial r_k}(r_k - \bar{r}_k)+ f_i(\bar{x},\bar{r}),~\forall i \in I
\end{align}
where $\frac{\partial f_i(\bar{x},\bar{r})}{\partial x_j}$ and  $\frac{\partial f_i(\bar{x},\bar{r})}{\partial r_k}$ are the partial derivative of $f_i$ with respect to $x_j$ and $r_k$ evaluated at $(\bar{x},\bar{r})$. Clearly, any feasible point satisfying~(\ref{eqt:beta}) is also feasible in the region defined by~(\ref{eqt:valid}), i.e., (\ref{eqt:valid}) does not eliminate any feasible region of P3. We can then model P3 using the following (sparse) MILP formulation:
\begin{align}
 \max~& \sum_{i \in I} 2 y_i \beta_i - \sum_{i \in I}y_i^2( \sum_{k \in K} \tilde{\theta}_{ik}r_k +  \sum_{j \in J} \theta_{ij} x_j ) \\
 st.~~&\beta_i  \leq  \sum_{j \in J}\frac{\partial f_i(x^t,r^t)}{\partial x_j} (x_j - x^t_j) +   \sum_{k \in K}\frac{\partial f_i(x^t,r^t)}{\partial r_k}(x_j - x^t_j) + f(x^t,r^t),~ \forall i \in I,  (x^{t},r^{t}) \in T \\
 \nonumber &(\ref{constr:P1-1}) - (\ref{constr:P1-4})
\end{align}
 where $T$ is the set of recorded points $(x,r)$. As the number of points in $T$ increases, the above formulation provides better approximation for P3.

With the above program, we now present our proposed algorithm in Algorithm~\ref{alg:QT_LA}. We refer to it as \textit{Quadratic Transform with Linear Alternating} (QT-LA) algorithm.

In Step 0, we initialize the solution and parameters. For simplicity, through this paper, we initialize $x^0_j = 1$,$\forall j \in J$, $r^0_k=1$,$\forall k \in K$, because the initial solution need not to be feasible. $n$ is the iteration number and $N_{max}$ is the maximum number of iterations.

In Step 1, we update set $T$ and $y$ and compute the approximated partial derivative. Here, $\gamma$ serves as the step-size. It is a number between 0 and 1, which controls the intensity of moving $y$ towards the value defined by  (\ref{eqt:opt_y}) at current solution $(x^n,y^n)$. In practice, this parameter is important for the algorithm. A brief discussion of its impacts is presented in~\ref{appendixB}. When computing the partial derivative, we add a small number $\epsilon$ to the denominate. This is to enhance the numerical stability and avoid the zero value in the denominate.

In Step 2, we solve a spare MILP.  This can be done using advanced MILP solvers. As the algorithm proceeds, new points (tentative solutions during iterations) are added to set $T$ and thus, the MILP contains more constraints and yields better approximation for the original formulation.

Finally, Step 3 is to check whether the stopping condition is met. The algorithm terminates when it finds repeated solutions or the maximum number of iterations exceeds a predetermined number. If the stopping condition is unmet, the process returns to Step 1.

The effectiveness and efficiency of QT-LA algorithm are discussed in Section~\ref{s:numeric_suggest}.

\begin{algorithm}[H] \label{alg:QT_LA}
\SetAlgoLined
Step 0: Initialize $x^0$ and $r^0$; $y^0 = 0$; $n=0$; $\epsilon = 0.0001$; $T = \emptyset$; choose $\gamma \in (0,1]$; $N_{max} \in Z_+$.\\
Step 1: $T: = T \cup \{(x^{n},r^{n})\}$.  $n := n +1$. Update parameters:
\begin{align*}
&y^n_i  := (1-\gamma)y^{n-1}_i + \gamma \cdot \frac{\sqrt[]{\sum_{k \in K} \tilde{b}_{ik}r^{n-1}_k + \sum_{j \in J} b_{ij}x^{n-1}_j}}{\sum_{k \in K}\tilde{\theta}_{ik}r^{n-1}_k +  \sum_{j \in J} \theta_{ij} x^{n-1}_j} , ~\forall i \in I \\
&\frac{\partial f^n_i}{\partial x_j} := \frac{b_{ij}}{2\sqrt[]{\sum_{k \in K} \tilde{b}_{ik}r^{n-1}_k + \sum_{j \in J} b_{ij}x^{n-1}_j} + \epsilon}, ~\forall i \in I, j \in J \\
&\frac{\partial f^n_i}{\partial r_k} := \frac{\tilde{b}_{ik}}{2\sqrt[]{\sum_{k \in K} \tilde{b}_{ik}r^{n-1}_k + \sum_{j \in J} b_{ij}x^{n-1}_j}+ \epsilon}, ~\forall i \in I, j \in J
\end{align*} \\
Step 2: Solve the following (sparse) MILP:
\begin{align*}
\max~&  \sum_{i \in I} 2 y^n_i \beta_i - \sum_{i \in I}(y_i^n)^2( \sum_{k \in K}\tilde{\theta}_{ik}r_k +  \sum_{j \in J} \theta_{ij} x_j ) \\
st.~& \beta_i  \leq  \sum_{j \in J}\frac{\partial f^t_i}{\partial x_j} (x_j - x^t_j) +   \sum_{k \in K}\frac{\partial f^t_i}{\partial r_k}(x_j - x^t_j) + f(x^t,r^t),~ \forall i \in I,  (x^{t},r^{t}) \in T \\
&(\ref{constr:P1-1}) - (\ref{constr:P1-4})
\end{align*}
to obtain the solution $(x^{n},r^{n})$. \\
Step 3: If $(x^{n},r^{n}) \in T$ or $n > N_{max}$, stop and output $(x^{n},r^{n})$. Else,  return to Step 1.
 \caption{QT-LA Algorithm}
\end{algorithm}

\subsection{Improve  procedure: particle swarm optimization}

In Improve procedure,  our goal is to improve the solution obtained from QT-LA. This can be done by using stochastic searching algorithms, e.g. genetic algorithm, variable neighborhood search, greedy randomized adaptive search procedure and particle swarm optimization (PSO). For P1, the decision variables are $x$ and $r$. The total number of dimensions is $|J|+|K|$. In practice, the above algorithms can be cast to handle problems with significantly larger variable dimensions. Therefore, the choice of the algorithm is \textit{flexible}. Meanwhile, there exist plentiful computer packages that enable straightforward and customized implementations of the algorithms. For instances, MATLAB has a global optimization toolbox that includes various metaheuristics for continuous and discrete problems. In Python, Pyevolve~\cite{perone2009pyevolve} and DEAP~\cite{fortin2012deap} both provide complete evolutionary algorithm frameworks. More recently, PySwarms~\cite{miranda2018pyswarms} enables basic and extendable optimization frameworks with PSO. These packages are powerful and can integrate with parallel computing technology to speed up the searching process. Therefore, we prefer to leverage them to design our algorithm rather than to write the whole program from scratch.

In this paper, we adopt the discrete binary version of PSO proposed by~\cite{kennedy1997discrete} and build our algorithm upon the PySwarms framework. We start by defining a list of binary variables $X = [r,X]$, i.e., $X$ is the concatenation of $r$ and $x$ with $X_k = r_k$,$\forall k \in K$ and $X_{j+|K|} = x_j$,$\forall j \in J$. Now, consider the unconstrained 0-1 problem with variable $X$:
\begin{align}
\textit{\textbf{[P4]}}~\label{eqt:PSO} \max~G(X) =\sum_{i \in I} \frac{ \sum_{k \in K} \tilde{b}_{ik}X_k +  \sum_{j \in J} b_{ij}X_{j +|K|}}{\sum_{k \in K}\tilde{\theta}_{ik}X_k +  \sum_{j \in J} \theta_{ij} X_{j +|K|}} - \varrho\left(\sum_{j \in J} X_{j+|K|} - P\right) - \varrho\left(|K| - P - \sum_{k \in K} X_k\right)
\end{align}
where $\varrho(y)$ is the penalty function. It can be chosen as $\varrho(y)=\max \{y,0\}$. We can verify that P4 is equivalent to P1. Essentially, we are heavily penalizing the violations of constraints (\ref{constr:P1-1}) - (\ref{constr:P1-2}) by the last two terms in (\ref{eqt:PSO}). When a point $\bar{X} = [\bar{r},\bar{x}]$ is not feasible in P1, the objective of P4 will be negative. In the optimal solution, both constraints will be satisfied and the objective of P4 will be equal to that of P1.

Without constraints, P4 can be solved using the standard binary PSO algorithm~\cite{kennedy1997discrete}. Here, $G(X)$, as in (\ref{eqt:PSO}), is the utility function to be maximized. Denote $S$ the number of particles in the swarm. Each particle has a position $X_s \in \{0,1\}^{|J|+|K|}$ and velocity $v_s \in \Re^{|J|+|K|}$. In our implementation, we set the initial positions of all particles to the solution obtained from QT-LA in Suggest procedure, while the velocity is initialized randomly. The impacts of running PSO after QT-LA are discussed in Section~\ref{s:numeric_suggest}.

\paragraph{Remark} We have presented the methodologies for Suggest-and-Improve framework. In Suggest procedure, we propose QT-LA algorithm. In Improve procedure, we apply PSO. From here onwards, we will refer to Suggest procedure as QT-LA and Suggest-and-Improve as QT-LA+PSO.

\section{Numerical study}
\label{s:ns}
In this section, we conduct numerical studies on the performance of different solution approaches. The results show that adding MC to MILP can largely reduce the computational time. However, for large-scale problems, MILP+MC fails to solve them in a resealable time, whereas our proposed QT-LA and QT-LA+PSO can efficiently handle them in seconds and leads to better solution quality (compared to the solution found by MILP+MC in 3600 seconds) in general. Finally, we present a case study based on the POP-Locker system in Singapore.

\subsection{Algorithm Performance}

We use artificially generated data. The Cartesian coordinates of the node set $I$, $J$ and $K$ are randomly generated by a uniform distribution in the interval $[0,1000]$. The demand of each node in $I$ are randomly generated by a uniform distribution in the interval $[0,100]$. The deterministic part of utility is defined as $v_{im} = -\alpha L_{im}$, where $L_{im}$ is the euclidian distance between node $i$ and $m$, $\forall i \in I, m \in J \cup K$. The unit of  $L_{im}$ is $100$. For the value of $a_{im}$, we adopt the following step-wise function: $a_{im} = 1$, if $L_{im} \leq 1$; $a_{im} = 0.5$, if $ 1<L_{im} \leq 2$; $a_{im} = 0.2$, if $2<L_{im} \leq 3$; $a_{im} = 0$, if $L_{im} > 3$.

We solve the MILPs and QT-LA through the Python-embedded modeling language CVXPY~\cite{cvxpy} and use Gurobi 8.1.1 under default setting as the solver. All computational experiments are done on a 16 GB memory iOS computer with 2.6 GHz Intel Core i7 processor.

\subsubsection{Performance of MILP formulations}
\label{s:numeric_mc}

We investigate the performance of MILP formulations and impact of adding MC. We use three artificial generated networks, namely, the 60-node network,  the 80-node network  and the 100-node network. The number of customer zones, POP-stations and candidate lockers is 30, 15 and 15 in the 60-node network; 40, 20 and 20 in the 80-node network; and 50, 25 and 25 in the 100-node network. We set the maximal computational time to 3600 seconds. Table~\ref{tab:milp_vs_mc} reports the results.

In general, the computational time declines with $\alpha$ increasing and grows with $P$ increasing. That is, for both formulations, when the value of $\alpha$ is small and $P$ is large, the problem instances require significantly more computational efforts than those with large $\alpha$ value and small $P$ value.

In almost all problem instances, adding the MC can lead to faster computation or smaller MIP gap, especially when the instances are difficult to solve by MILP.  For example, in the 80-node network, when $\alpha = 0.5$ and $P=5$, MILP requires 1451.8s to find the global optimal solution, while MILP+MC requires only 114.5s. 
Meanwhile, when the instances cannot be solved, MILP+MC terminates with much smaller gaps. The better performance of MILP+MC is due to the stronger continuous relaxation. Therefore, when implementing the MILP approach, it is important to add the MC. 

Finally, the results suggest both formulations cannot scale well on the problem size. In the 60-node network, all problem instances can be solved easily using both formulations. The maximum computational times among these 18 instances are 150.6s and 55s for MILP and MILP+MC. However, when we slightly increase the problem size to 80 nodes (add 10 customer zones, 5 pos-stations and 5 lockers), MILP fails to solve the model within one hour in some problem instances. The computational time for MILP+MC also grows dramatically.  In the 100-node network, when $\alpha = 0.5$ and $\alpha = 1$, both formulations are not efficiently in solving these instances. Therefore, MILP approaches are not suitable for large-scale problems.

\begin{table}[]\caption{Computational time and Gap via MILP and MILP+MC under the 60-node, 80-node and 100-node problem instances. CPU denotes the times in seconds used to solve the respective problem. The maximal computational time is set to 3600 seconds. $**$ indicates the instance cannot be solved to global optimality. Gap is the MIP gap reported by Gurobi in \%.}  \label{tab:milp_vs_mc}
\centering
\scriptsize
\begin{tabular}{lllllllllllllllllll}
\hline
\multicolumn{1}{c}{\multirow{3}{*}{$\alpha$}} & \multirow{3}{*}{P} & \multicolumn{5}{c}{60-node}                              &  & \multicolumn{5}{c}{80-node}                              &  & \multicolumn{5}{c}{100-node}                             \\ \cline{3-7} \cline{9-13} \cline{15-19}
\multicolumn{1}{c}{}                   &                    & \multicolumn{2}{c}{MILP} &  & \multicolumn{2}{c}{MILP+MC} &  & \multicolumn{2}{c}{MILP} &  & \multicolumn{2}{c}{MILP+MC} &  & \multicolumn{2}{c}{MILP} &  & \multicolumn{2}{c}{MILP+MC} \\ \cline{3-4} \cline{6-7} \cline{9-10} \cline{12-13} \cline{15-16} \cline{18-19}
\multicolumn{1}{c}{}                   &                    & CPU          & Gap       &  & CPU           & Gap         &  & CPU         & Gap        &  & CPU            & Gap        &  & CPU         & Gap        &  & CPU           & Gap         \\ \hline
0.5    & 1                  & 0.2          & 0         &  & 0.3           & 0           &  & 0.4         & 0          &  & 0.6            & 0          &  & 0.6         & 0          &  & 0.7           & 0           \\
         & 3                  & 2.9          & 0         &  & 1.3           & 0           &  & 15.4        & 0          &  & 3.4            & 0          &  & 10.0        & 0          &  & 4.1           & 0           \\
         & 5                  & 52.8         & 0         &  & 8.8           & 0           &  & 1451.8      & 0          &  & 114.5          & 0          &  & $**$      & 0.13       &  & 264.5         & 0           \\
         & 7                  & 115.9        & 0         &  & 38.3          & 0           &  & $**$      & 9.88       &  & 592.1          & 0          &  & $**$      & 16.42      &  & $**$        & 2.68        \\
         & 9                  & 150.6        & 0         &  & 55.0          & 0           &  & $**$      & 12.22      &  & 1460.6         & 0          &  & $**$      & 25.06   &  & $**$        & 7.96        \\
         & 10                 & 93.8         & 0         &  & 43.2          & 0           &  & $**$      & 13.35      &  & 1908.3         & 0          &  & $**$      & 27.61   &  & $**$        & 9.91        \\
\\
1       & 1                  & 0.4          & 0         &  & 0.2           & 0           &  & 0.9         & 0          &  & 0.7            & 0          &  & 0.7         & 0          &  & 0.7           & 0           \\
         & 3                  & 3.5          & 0         &  & 1.2           & 0           &  & 20.9        & 0          &  & 6.8            & 0          &  & 29.1        & 0          &  & 7.1           & 0           \\
         & 5                  & 17.6         & 0         &  & 2.7           & 0           &  & 371.8       & 0          &  & 86.3           & 0          &  & 3506.1      & 0          &  & 734.4         & 0           \\
         & 7                  & 19.7         & 0         &  & 3.8           & 0           &  & 1833.9      & 0          &  & 227.3          & 0          &  & $**$      & 9.70       &  & 2683.7        & 0           \\
         & 9                  & 14.0         & 0         &  & 3.1           & 0           &  & 2696.3      & 0          &  & 272.7          & 0          &  & $**$      & 11.09      &  & $**$        & 0.11        \\
         & 10                 & 13.4         & 0         &  & 3.2           & 0           &  & 2351.9      & 0          &  & 242.9          & 0          &  & $**$     & 10.21      &  & 2940.7        & 0           \\
\\
2      & 1                  & 0.2          & 0         &  & 0.2           & 0           &  & 0.8         & 0          &  & 0.5            & 0          &  & 1.0         & 0          &  & 0.5           & 0           \\
        & 3                  & 2.3          & 0         &  & 1.8           & 0           &  & 8.7         & 0          &  & 4.6            & 0          &  & 7.7         & 0          &  & 8.7           & 0           \\
        & 5                  & 1.7          & 0         &  & 1.7           & 0           &  & 16.5        & 0          &  & 23.1           & 0          &  & 64.2        & 0          &  & 77.5          & 0           \\
        & 7                  & 2.0          & 0         &  & 1.7           & 0           &  & 197.7       & 0          &  & 33.1           & 0          &  & 1467.0      & 0          &  & 186.8         & 0           \\
        & 9                  & 1.4          & 0         &  & 1.0           & 0           &  & 312.7       & 0          &  & 69.5           & 0          &  & 2151.3      & 0          &  & 151.3         & 0           \\
        & 10                 & 1.9          & 0         &  & 0.8           & 0           &  & 315.7       & 0          &  & 119.4          & 0          &  & 1808.3      & 0          &  & 146.6         & 0           \\ \hline
\end{tabular}
\end{table}

\subsubsection{Performance of Suggest-and-Improve framework}
\label{s:numeric_suggest}

We test the performance of our proposed QT-LA and QT-LA+PSO. We use two artificial generated networks, namely, the 100-node (small) network and the 400-node  (large) network. The number of customer zones, POP-stations and candidate lockers is 50, 25 and 25 in the 100-node network; and 200, 100 and 100 in the 400-node network.

We first compare the performance of QT-LA with MILP+MC. We set the maximal computational time for MILP+MC to 3600s. Table~\ref{tab:suggest_improve} reports the computational results of 24 problem instances. For QT-LA, $\gamma^*$ denotes the best step size, which is selected from $\{0.4,0.5,0.6,0.7,0.8,0.9,1\}$. $\Delta^1$ denotes the relative difference between the objective by MILP+MC and the objective by QT-LA. It is computed by $(C-C^1)/C^1$. A positive value of $\Delta^1$ indicates MILP+MC finds a better solution than QT-LA; vice versa.

For small-scale problems, QT-LA is much faster and yields high-quality solutions. In the 100-node network, MILP+MC is able to solve most of the instances to global optimality. In this case, the maximum $\Delta^1$ is only $0.29\%$, meaning that QT-LA finds near-optimal solutions with minimal computational efforts (CPU is less than 0.5s).  For large-scale problems, QT-LA is a more reliable method that requires significantly less computation load and leads to better solution quality as well. In the 400-node network, $\Delta^1$ is negative in all instances. In the worst case, $\Delta^1$ drops to $-10.65\%$. The computational time by QT-LA is still not significant.


\begin{table}[]\caption{Performance of MILP+MC, QT-LA and QT-LA+PSO. $N$ denotes the number of nodes. $C$, $C^1$, $C_{avg}$ and $C_{max}$ denote the service level (\%). CPU denotes the time in seconds used to solve the problem. $**$ indicates that the instance cannot be solved within 3600s by MILP+MC.}  \label{tab:suggest_improve}
 \centering \scriptsize
\begin{tabular}{lllllllllllllll}
\hline
\multirow{2}{*}{N} & \multicolumn{1}{c}{\multirow{2}{*}{$\alpha$}} & \multirow{2}{*}{P} & \multicolumn{2}{c}{MILP+MC} & \multicolumn{1}{c}{} & \multicolumn{4}{c}{QT-LA}                   &  & \multicolumn{4}{c}{QT-LA+PSO}                  \\ \cline{4-5} \cline{7-10} \cline{12-15} 
                   & \multicolumn{1}{c}{}                          &                    & $C$          & CPU          &                      & $C^1$ & $\gamma^*$ & CPU   & $\Delta^1$(\%) &  & $C_{avg}$ & $C_{max}$ & CPU   & $\Delta^2$(\%) \\ \cline{1-10} \cline{12-15} 
100                & 0.5                                           & 5                  & 31.85        & 264.5        &                      & 31.85 & 0.8        & 0.1   & 0              &  & 31.85     & 31.85     & 4.9   & 0              \\
                   &                                               & 10                 & 34.58        & $**$         &                      & 34.58 & 0.8        & 0.1   & 0              &  & 34.58     & 34.58     & 5.1   & 0              \\
                   &                                               & 15                 & 36.66        & $**$         &                      & 36.76 & 0.8        & 0.1   & -0.27          &  & 36.76     & 36.76     & 5.1   & -0.27          \\
                   &                                               & 20                 & 37.63        & $**$         &                      & 37.63 & 0.4        & 0.5   & 0              &  & 37.75     & 37.76     & 5.3   & -0.64          \\
                   & 1                                             & 5                  & 52.24        & 734.4        &                      & 52.21 & 0.8        & 0.1   & 0.06           &  & 52.22     & 52.24     & 5.1   & 0              \\
                   &                                               & 10                 & 55.42        & 2940.7       &                      & 55.42 & 0.8        & 0.2   & 0              &  & 55.42     & 55.42     & 5.2   & 0              \\
                   &                                               & 15                 & 58.14        & 780.0        &                      & 58.14 & 0.8        & 0.1   & 0              &  & 58.14     & 58.14     & 5.0   & 0              \\
                   &                                               & 20                 & 58.31        & 1242.1       &                      & 58.31 & 0.6        & 0.3   & 0              &  & 58.31     & 58.31     & 5.2   & 0              \\
                   & 2                                             & 5                  & 71.69        & 77.5         &                      & 71.66 & 0.8        & 0.1   & 0.04           &  & 71.67     & 71.69     & 5.1   & 0              \\
                   &                                               & 10                 & 73.7         & 146.6        &                      & 73.65 & 0.8        & 0.2   & 0.07           &  & 73.7      & 73.7      & 5.1   & 0              \\
                   &                                               & 15                 & 74.67        & 171.8        &                      & 74.55 & 0.5        & 0.3   & 0.16           &  & 74.67     & 74.67     & 5.3   & 0              \\
                   &                                               & 20                 & 74.69        & 1342.1       &                      & 74.47 & 0.5        & 0.3   & 0.29           &  & 74.69     & 74.69     & 5.3   & 0              \\
                   &                                               &                    &              &              &                      &       &            &       &                &  &           &           &       &                \\
400                & 0.5                                           & 10                 & 30.69        & $**$         &                      & 31.83 & 0.7        & 1.2   & -3.38          &  & 31.83     & 31.83     & 11.9  & -3.38          \\
                   &                                               & 20                 & 29.54        & $**$         &                      & 33.06 & 0.7        & 1.4   & -10.65         &  & 33.06     & 33.06     & 12.3  & -10.65         \\
                   &                                               & 40                 & 31.39        & $**$         &                      & 34.19 & 0.9        & 2.5   & -8.19          &  & 34.19     & 34.19     & 13.3  & -8.19          \\
                   &                                               & 60                 & 33.48        & $**$         &                      & 35.05 & 0.7        & 5.1   & -4.48          &  & 35.05     & 35.05     & 16.4  & -4.48          \\
                   & 1                                             & 10                 & 54.17        & $**$         &                      & 54.24 & 0.9        & 1.1   & -0.13          &  & 54.24     & 54.24     & 11.6  & -0.13          \\
                   &                                               & 20                 & 51.17        & $**$         &                      & 55.93 & 0.9        & 1.2   & -8.51          &  & 55.93     & 55.93     & 12.2  & -8.51          \\
                   &                                               & 40                 & 52.96        & $**$         &                      & 57.37 & 0.9        & 8.5   & -7.69          &  & 57.37     & 57.37     & 19.8  & -7.69          \\
                   &                                               & 60                 & 55.45        & $**$         &                      & 58.38 & 1          & 63.4  & -5.02          &  & 58.38     & 58.39     & 77.7  & -5.04          \\
                   & 2                                             & 10                 & 76.61        & $**$         &                      & 78.99 & 0.8        & 1.6   & -3.01          &  & 78.99     & 79.0      & 12.2  & -3.02          \\
                   &                                               & 20                 & 76.77        & $**$         &                      & 80.43 & 0.7        & 7.6   & -4.55          &  & 80.43     & 80.44     & 18.2  & -4.58          \\
                   &                                               & 40                 & 77.85        & $**$         &                      & 81.72 & 0.9        & 44.8 & -4.74          &  & 81.73     & 81.76     & 58.0  & -4.78          \\
                   &                                               & 60                 & 79.54        & $**$         &                      & 82.31 & 1          & 107.7 & -3.37          &  & 82.52     & 82.59     & 121.5 & -3.69          \\ \hline
\end{tabular}
\end{table}

Next, we discuss the impacts of adding the Improve procedure. For all instances, we implemented QT-LA+PSO algorithm for $10$ replications. In PSO, the number of particles and iterations in each replication are $50$ and $2000$. The results are shown in Table~\ref{tab:suggest_improve}. Here, CPU is the average computational time. $C_{avg}$ is the average objective and $C_{max}$ is the best objective. $\Delta^2$ denotes the relative difference between the objective by MILP+MC and the best objective by QT-LA+PSO. It is computed by  $(C-C_{max})/C_{max}$. A positive value of $\Delta^2$ indicates MILP+MC finds a better solution than QT-LA+PSO; vice versa.

From Table~\ref{tab:suggest_improve}, running the PSO after QT-LA will only slightly increase the computational load. In terms of solution quality, it provides a chance to further improve the solution. We observe that, in all instances, the solutions by QT-LA+PSO are no worse than the solutions by MILP+MC ($\Delta^2 \leq 0$).

To conclude, QT-LA+PSO is an efficient algorithm that is suitable for large-scale problems.

\subsection{Case Study}

To illustrate the usefulness of our model, we conduct a case study based on the POP-station network and demand zones of Singapore.  

\subsubsection{Data collection and the impact of $P$}

Four data sets are described below.
\begin{itemize}
\item Demand and demand zone. We assume the demand comes from the residential area. Demand nodes were formed based on the sub-clusters demarcated by Singapore’s Urban Redevelopment Authority (URA). The size of each cluster and the number of residents are retrieved from a publication by the Department of Statistics (DOS) Singapore~\footnote{https://www.singstat.gov.sg/-/media/files/publications/population/population2016.pdf}. There are 272 residential clusters. The demand is assumed proportional to the number of residents.
\item POP-stations. 137 POP-stations located around Singapore are identified. The address of each POP-station is obtained from Singapore Post's official webpage.
\item Potential locker facilities. 139 shopping centers are chosen as potential locations. The list of shopping centers is extracted from Shopping.sg, a shopping guide that provides extensive information on shopping centers in Singapore. Using Google maps, the addresses of these shopping centers are  retrieved and added to the list.
\item Distance matrix $L$. We use the Google's API URL to calculate accurate walking distances between nodes. The distance is measured in kilometers with 3 decimal.
\end{itemize}
That is, $|I| = 272$, $|J| = 139$ and $|K| = 137$. The deterministic part of utility is defined as $v_{im} = -\alpha L_{im}$, $\forall i \in I, m \in J \cup K$. 
For the value of $\tilde{a}$, $\tilde{a}_{im} = 1$, if $L_{im} \leq 1$; $\tilde{a}_{im} = 0.5$, if $ 1<L_{im} \leq 1.5$; $\tilde{a}_{im} = 0.2$, if $1.5<L_{im} \leq 2$; and $\tilde{a}_{im} = 0$, if $L_{im} > 2$.


We first consider the problem where $P$ POP-stations shall be closed and $P$ locker facilities are to be located, namely,  (\ref{constr:P1-1}) and (\ref{constr:P1-2}) are equality constraints. In this case, we can use $\varrho(y)= y^2$ as the penalty function in PSO. We test the model under different $P$ values and $\alpha$ values.  For each problem instances, we implement QT-LA+PSO with tuned $\gamma$ for 20 replications and selected the best solution from the pool. For PSO, the number of particles and iterations in each replication are $50$ and $5000$.

\begin{figure}
  \centering
    \includegraphics[width=80mm]{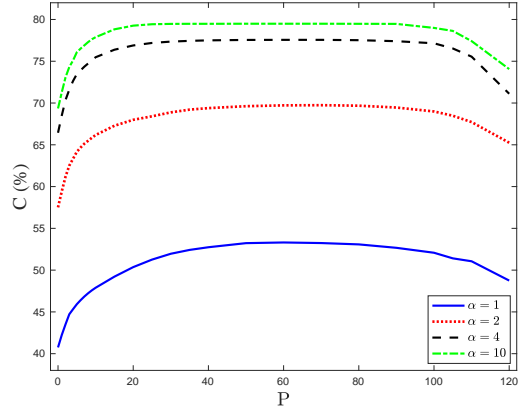}
      \caption{Case study. Service level $C$ versus $P$.} \label{fig:case_study_result}
\end{figure}

Figure~\ref{fig:case_study_result} depicts the results. The service level $C$ first increases with $P$ to a certain maximum and then decreases. The maximum improvement in the service level is around $10\%-12\%$, which occurs when $P$ is near $60$ for all values of $\alpha$. This implies that some candidate lockers are located in the sites closer to the customer zones than the POP-stations. Opening these lockers instead of some remote POP-stations will improve the overall service level.

In general, the service level hardly improves when $P$ increases from 25 to 60. For examples, under $\alpha=4$, $C$ values for $P=25$ and $P=60$ are $77.19\%$ versus $77.56 \%$. Under $\alpha =10$, $C$ values for $P=25$ and $P=60$ are $79.43\%$ versus $79.49 \%$. These results suggest that replacing 60 existing POP-stations with 60 carefully selected locker facilities will not lead to an acceptable increase in the service level compared to only replacing 25 POP-stations. In practice, the latter option is advantageous  because the effort for making the change is much smaller.  On the other hand, customers may be accustomed to the current network and facilities. Closing a large number of POP-stations may result in their discomfort because their frequently visited POP-stations may be closed and they have to adjust the change. Therefore, it is reasonable to replace only a moderate number of POP-stations. From Figure~\ref{fig:case_study_result}, an appropriate $P$ value can be selected from $20$ to $25$.

\subsubsection{Impact of $\alpha$}

Next, we discuss the impact of $\alpha$ on the service level and the solution.  From Figure~\ref{fig:case_study_result}, when $\alpha$ increases, $C$ will increase. Intuitively, a large value of $\alpha$ implies that the customers show high preference to their nearest facility, meaning that only a small proportion of customers will visit the non-nearest facilities. The service level provided by the alliance will thus be high.

Let $\tilde{K}_{\alpha}=\{k \mid r_k=0, \forall k \in K\}$, that is, $\tilde{K}_{\alpha}$ is the set of closed POP-stations suggested by the model. Table~\ref{tab:closed_pop} shows the sets corresponding to 3 different values of $P$. We see that the solutions are different when $\alpha$ values are different. In particular, for all values of $P$, $\tilde{K}_1 \cap \tilde{K}_{10} = \emptyset$, indicating that $\tilde{K}_1$ is completely different from $\tilde{K}_{10}$. Therefore, in order to make the best decision in practice, the alliance needs to carefully estimate the parameter $\alpha$.

\begin{table}[]\caption{The set of closed POP-stations suggested by the solution, namely, $\tilde{K}_{\alpha}$.}  \label{tab:closed_pop}
 \centering \scriptsize
\begin{tabular}{llllllll}
\hline
P & $\alpha=1$              &  & $\alpha=2$            &  & $\alpha=4$              &  & $\alpha=10$             \\ \hline
1 & \{62\}                  &  & \{88\}                &  & \{136\}                 &  & \{132\}                 \\
3 & \{15, 33, 88\}          &  & \{2, 88, 136\}        &  & \{41, 93, 136\}         &  & \{93, 132, 135\}        \\
5 & \{15, 33, 82, 88, 129\} &  & \{2, 9, 36, 88, 136\} &  & \{9, 41, 93, 131, 136\} &  & \{0, 9, 120, 135, 136\} \\ \hline
\end{tabular}
\end{table}

In our final experiment, we consider the scenario where the alliance wants to close 20 POP-stations and open an unrestricted number of locker facilities to maximize the service level, that is, we set $P=20$ and remove constraint~(\ref{constr:P1-1}). As shown in Table~\ref{tab:p20}, as $\alpha$ increases, the model suggests to open more lockers and the total number of facilities increases. This is reasonable because the customers' choice behaviors are ``random" in the case of a small $\alpha$ value.  Given too many facilities, customers will ``randomly" select one to visit.  In many cases, they may walk a long distance, resulting in a lower service level. Conversely, with a large $\alpha$ value, customers will almost always patronize the nearest facility. It is then reasonable to operate more facilities (than the case with a small $\alpha$ value) to gain proximity to customer zones.

According to Table~\ref{tab:p20}, the model does not recommend to open $|J|$ lockers. This result is fundamentally different from the traditional coverage model where the solution is to open $|J|$ lockers if the maximum number of lockers (that are allowed to open) is unrestricted and the cost is not considered. The reason for the difference lies in the model assumptions. More specifically, the traditional coverage model assumes that the operator determines the customer-facility allocation or that customers always select the nearest facility. However, this assumption is not valid in e-commerce because customers take the initiative. Meanwhile, the choice behaviors are typically probabilistic rather than deterministic~\cite{gul2014random,Haase2014comparison,ljubic2018outer}. Therefore, the MNL based on random utility theory is applied to forecast the choice. Under MNL, if we open an excess number of facilities, the chance of a customer patranizing the nearest one will decreases. In some cases, we could observe an interesting phenomenon, that is, the alliance intends to enhance the service by opening new facilities, but the result is the opposite. 

\begin{table}[]\caption{Results under $P=20$ when constraint~(\ref{constr:P1-1}) is removed. $C$ denotes the service level (\%); $X$ denotes the number of lockers; and $F$ denotes the total number of facilities.}  \label{tab:p20}
 \centering \scriptsize
\begin{tabular}{llllllll}
\hline
      & $\alpha=1$ &  & $\alpha=2$ &  & $\alpha=4$ &  & $\alpha=10$ \\ \hline
$C$   & 50.54      &  & 68.07      &  & 77.18      &  & 79.48       \\
$X$ & 29         &  & 40         &  & 44         &  & 51          \\
$F$   & 146        &  & 157        &  & 161        &  & 168         \\ \hline
\end{tabular}
\end{table}

To summarize, under different $\alpha$ values, the final decisions could be completely different. When determining where to locate the locker facilities and which POP-stations to close, it is important to consider the customers' choices. The value of $\alpha$, therefore, should be carefully determined or estimated in advance (e.g. conducting a parameter estimation using survey data). Otherwise, we may operate facilities that ultimately reduce the overall service level.

\section{Conclusion}
\label{s:conclusion}

This paper focused on the last-mile delivery problem. We considered a POP-Locker Alliance who operates a set of POP-stations and wishes to improve the service level by opening new locker facilities. We studied the optimal location of locker facilities, considering the customer' preference, with the objective to maximize the overall service level. The customer's choice was modeled by a MNL model. We formulated the problem as a MLFP. To solve the model, we provided two approaches. The first one was the MILP approach. We reformulated our original problem as an equivalent MILP formulation and further strengthened the formulation using McCormick inequalities (MC). Our numerical test suggested that adding MC can significantly reduce the computational time. However, the strengthened MILP cannot scale well on the problem size. To solve large-scale problems, we applied the Suggest-and-Improve framework. In the Suggest procedure, we found a high-quality solution using our proposed QT-LA algorithm. The solution by QT-LA was then improved by a PSO algorithm in Improve procedure. In our numerical tests, our framework outperformed MILP approach in both solution time and solution quality. Finally, we conducted a case study to illustrate the usefulness of our model. The results highlighted the importance of considering the customers' choice behaviors and the difference between traditional coverage model and our model. 

There are limitations. We did not consider the facility capacity. Indeed, capacity is an important factor that can affect the customer's choice behavior. For example, if a locker has no remaining capacity, customers may select other facilities. We assumed that customers are homogeneous in their observable characteristics. This strong assumption can be partially relaxed by dividing customers into different classes according to their characteristics such as age, income and car-ownership. Besides, we measure the service level by a function of distance and the customer's choices are developed based on distance as well. In practice, other factors, such as transportation convenience and proximity to working places, are also important. Future research can consider them in the choice model and propose a more general measurement for the service level.  Finally, for large-scale problems, we applied QT-LA+PSO. Even though it works well in practice, it is not an exact solution approach. Future studies can develop an exact and efficient solution approach.

%
%
\appendix
\section{}
\label{appendix}

\subsection{Proof of Proposition 1}
\begin{proof}
The upper bound on $z_i$ is equivalent to the lower bound on ${\sum_{k \in K}\tilde{\theta}_{ik}r_k +  \sum_{j \in J} \theta_{ij} x_j}$. For $i \in I$, we define the auxiliary problem:
\begin{align*}
h_i = &\min~ \sum_{k \in K}\tilde{\theta}_{ik}r_k +  \sum_{j \in J} \theta_{ij} x_j \\
        st.~& (\ref{constr:P1-1}) - (\ref{constr:P1-4})
 \end{align*}
By (\ref{constr:P1-1}), $h_i$ is attained only when $x_j = 0, \forall j$. By (\ref{constr:P1-2}), we need to select at least $|K|-P$ elements from $\tilde{\theta}_{in}$,$\forall n \in K$. To attain $h_i$, we select the smallest $|K|-P$ elements.
\end{proof}

\subsection{Proof of Proposition 2}
\begin{proof}
For $\forall i \in I, j\in J$, define the following auxiliary problem:
\begin{align*}
h_{i|x_j = \xi} =~&\text{optimize}~ \sum_{k \in K}\tilde{\theta}_{ik}r_k +  \sum_{j \in J} \theta_{ij} x_j \\
        st.~& (\ref{constr:P1-1}) - (\ref{constr:P1-4}) \\
              & x_j = \xi
\end{align*}

For $z^u_{i \mid x_j=1}$, we set ``optimize" to ``min" and $\xi = 1$. The optimal value of this problem is attained only when $x_n = 0$,$\forall n \in J\backslash \{j\}$. By (\ref{constr:P1-2}), we need to select at least $|K|-P$ elements from $\tilde{\theta}_{in}$,$\forall n \in K$. To attain $h_i$, we select the smallest $|K|-P$ elements.

For $z^l_{i \mid x_j=0}$, we set  ``optimize" to ``max" and $\xi = 0$. The optimal value is attained only when $r_k = 1$,$\forall k \in K$. By (\ref{constr:P1-1}) and $x_j = 0$, we need to select at most $P$ elements from $\theta_{im}$,$\forall m \in J\backslash \{j\}$, to maximize the objective function, i.e, we should select the largest $P$ elements.

Similarly, for $z^l_{i \mid x_j=1}$,  we set  ``optimize" to ``max" and $\xi = 1$. The optimal value is attained only when $r_k = 1$,$\forall k \in K$. By (\ref{constr:P1-1}) and $x_j = 1$, we need to select the largest $P-1$ elements from $\theta_{im}$,$\forall m \in J\backslash \{j\}$.
\end{proof}

\subsection{Proof of Proposition 3}
\begin{proof}
For $\forall i \in I, k\in K$, define the following auxiliary problem:
\begin{align*}
h_{i|r_k = \xi} =~&\text{optimize}~ \sum_{k \in K}\tilde{\theta}_{ik}r_k +  \sum_{j \in J} \theta_{ij} x_j \\
        st.~& (\ref{constr:P1-1}) - (\ref{constr:P1-4}) \\
              & r_k = \xi
\end{align*}

For $z^u_{i \mid r_k=1}$, we set ``optimize" to ``min" and $\xi = 1$. The optimal value of this problem is attained only when $x_n = 0$,$\forall n \in J$. By (\ref{constr:P1-2}) and $r_k = 1$, we need to select at least $|K|-P-1$ elements from $\tilde{\theta}_{in}$,$\forall n \in K\backslash \{k\}$. To attain $h_i$, we select the smallest $|K|-P-1$ elements.

For $z^l_{i \mid r_k=0}$, we set  ``optimize" to ``max" and $\xi = 0$. The optimal value is attained only when $r_k = 1$,$\forall k \in K\backslash \{k\}$. By (\ref{constr:P1-1}), we need to select at most $P$ elements from $\theta_{im}$,$\forall m \in J$, i.e, to select the largest $P$ elements.

Similarly, for $z^l_{i \mid r_k=1}$,  we set  ``optimize" to ``max" and $\xi = 1$. The optimal value is attained only when we set $r_k = 1$,$\forall k \in K$, and select the largest $P$ elements from $\theta_{im}$,$\forall m \in J$.
\end{proof}

\section{}
\label{appendixB}
\begin{figure*}[h]
\begin{center}
     \psfig{figure=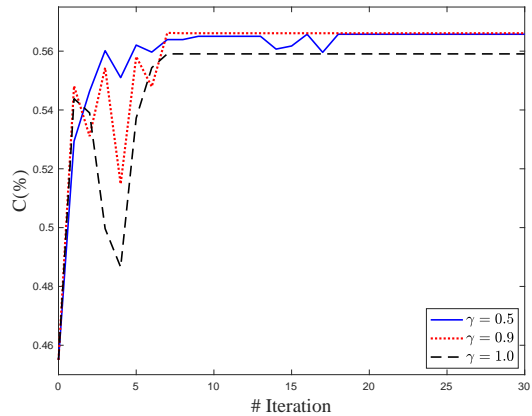,width=80mm}
\end{center}
\caption{Convergence of QT-LA under different $\gamma$ values. \# Iteration denotes the number of iterations.} \label{fig:learning_rate}
\end{figure*}

We briefly discuss the impact of $\gamma$ on QT-LA algorithm. The test problem instance is a 140-node network with $\alpha=1$ and $P=30$. In terms of the objective function value, the optimal $\gamma^*$ is $0.9$. QT-LA algorithm is run under 3 values of $\gamma$, namely, $0.5$ (low), $0.9$ (optimal), and $1.0$(high).  Figure~\ref{fig:learning_rate} shows the convergence of QT-LA algorithm.

For a smaller $\gamma$, the algorithm proceeds more smoothly. Intuitively, given a small step size, the changes in parameters will not be dramatical in two consecutive iterations. This also implies the algorithm may need more iterations to converge. From Figure~\ref{fig:learning_rate}, QT-LA converges in iteration 8 when $\gamma = 1$, in iteration 8 when $\gamma = 0.9$ and in iteration 18 when $\gamma = 0.5$.

Meanwhile, under different $\gamma$ values, QT-LA may converge to different solutions. It is thus important to tune the value of  $\gamma$. One method is to try different values and select the one that returns the best objective. From our numerical studies, this simple idea works well.

\bibliographystyle{plain}
\bibliography{references}
\end{document}